\documentclass[12pt]{article}
\usepackage{latexsym, amssymb}

\DeclareMathAlphabet\gothic{U}{euf}{m}{n}

\makeatletter
\def\eqnarray{\stepcounter{equation}\let\@currentlabel=\theequation
\global\@eqnswtrue
\tabskip\@centering\let\\=\@eqncr
$$\halign to \displaywidth\bgroup\hfil\global\@eqcnt\z@
  $\displaystyle\tabskip\z@{##}$&\global\@eqcnt\@ne
  \hfil$\displaystyle{{}##{}}$\hfil
  &\global\@eqcnt\tw@ $\displaystyle{##}$\hfil
  \tabskip\@centering&\llap{##}\tabskip\z@\cr}

\def\endeqnarray{\@@eqncr\egroup
      \global\advance\c@equation\m@ne$$\global\@ignoretrue}

\def\@yeqncr{\@ifnextchar [{\@xeqncr}{\@xeqncr[5pt]}}
\makeatother

\begin{document}
\bibliographystyle{tom}

\newtheorem{lemma}{Lemma}[section]
\newtheorem{thm}[lemma]{Theorem}
\newtheorem{cor}[lemma]{Corollary}
\newtheorem{voorb}[lemma]{Example}
\newtheorem{rem}[lemma]{Remark}
\newtheorem{prop}[lemma]{Proposition}
\newtheorem{stat}[lemma]{{\hspace{-5pt}}}
\newtheorem{obs}[lemma]{Observation}
\newtheorem{defin}[lemma]{Definition}

\newenvironment{remarkn}{\begin{rem} \rm}{\end{rem}}
\newenvironment{exam}{\begin{voorb} \rm}{\end{voorb}}
\newenvironment{defn}{\begin{defin} \rm}{\end{defin}}
\newenvironment{obsn}{\begin{obs} \rm}{\end{obs}}

\newenvironment{emphit}{\begin{itemize} }{\end{itemize}}

\newcommand{\gota}{\gothic{a}}
\newcommand{\gotb}{\gothic{b}}
\newcommand{\gotc}{\gothic{c}}
\newcommand{\gote}{\gothic{e}}
\newcommand{\gotf}{\gothic{f}}
\newcommand{\gotg}{\gothic{g}}
\newcommand{\gothh}{\gothic{h}}
\newcommand{\gotk}{\gothic{k}}
\newcommand{\gotm}{\gothic{m}}
\newcommand{\gotn}{\gothic{n}}
\newcommand{\gotp}{\gothic{p}}
\newcommand{\gotq}{\gothic{q}}
\newcommand{\gotr}{\gothic{r}}
\newcommand{\gots}{\gothic{s}}
\newcommand{\gotu}{\gothic{u}}
\newcommand{\gotv}{\gothic{v}}
\newcommand{\gotw}{\gothic{w}}
\newcommand{\gotz}{\gothic{z}}
\newcommand{\gotA}{\gothic{A}}
\newcommand{\gotB}{\gothic{B}}
\newcommand{\gotG}{\gothic{G}}
\newcommand{\gotL}{\gothic{L}}
\newcommand{\gotS}{\gothic{S}}
\newcommand{\gotT}{\gothic{T}}

\newcounter{teller}
\renewcommand{\theteller}{\Roman{teller}}
\newenvironment{tabel}{\begin{list}%
{\rm \bf \Roman{teller}.\hfill}{\usecounter{teller} \leftmargin=1.1cm
\labelwidth=1.1cm \labelsep=0cm \parsep=0cm}
                      }{\end{list}}

\newcounter{tellerr}
\renewcommand{\thetellerr}{(\roman{tellerr})}
\newenvironment{subtabel}{\begin{list}%
{\rm  (\roman{tellerr})\hfill}{\usecounter{tellerr} \leftmargin=1.1cm
\labelwidth=1.1cm \labelsep=0cm \parsep=0cm}
                         }{\end{list}}
\newenvironment{ssubtabel}{\begin{list}%
{\rm  (\roman{tellerr})\hfill}{\usecounter{tellerr} \leftmargin=1.1cm
\labelwidth=1.1cm \labelsep=0cm \parsep=0cm \topsep=1.5mm}
                         }{\end{list}}

%%%%%%%%%%%%%%%%

%\newcounter{teller}
%\renewcommand{\theteller}{\Roman{teller}}
%\newenvironment{tabel}{\begin{list}%
%{\rm \bf \Roman{teller}.\hfill}{\usecounter{teller} \leftmargin=1.1cm
%\labelwidth=1.1cm \labelsep=0cm \parsep=0cm}
%                     }{\end{list}}

\newenvironment{tabelpairs}{\begin{list}%
{{\rm \bf \Roman{teller}.}\ {\rm( {\bf \Roman{teller}$^\prime.\,$})}\hfill}{\usecounter{teller} \leftmargin=1.9cm
\labelwidth=1.9cm \labelsep=0cm \parsep=0cm}
                     }{\end{list}}

%\newenvironment{tabelpairs}{\begin{list}%
%{\rm \bf \Roman{teller}.\ (\Roman{teller}$^\prime$).\hfill}{\usecounter{teller} \leftmargin=1.9cm
%\labelwidth=1.9cm \labelsep=0cm \parsep=0cm}
%                     }{\end{list}}

%%%%%%%%%%%%%%

\newcommand{\Ni}{{\bf N}}
\newcommand{\Ri}{{\bf R}}
\newcommand{\Ci}{{\bf C}}
\newcommand{\Ti}{{\bf T}}
\newcommand{\Zi}{{\bf Z}}
\newcommand{\Fi}{{\bf F}}

\newcommand{\proof}{\mbox{\bf Proof} \hspace{5pt}} 
\newcommand{\remark}{\mbox{\bf Remark} \hspace{5pt}}
\newcommand{\ruimte}{\vskip10.0pt plus 4.0pt minus 6.0pt}

\newcommand{\simh}{{\stackrel{{\rm cap}}{\sim}}}
\newcommand{\ad}{{\mathop{\rm ad}}}
\newcommand{\Ad}{{\mathop{\rm Ad}}}
\newcommand{\Aut}{\mathop{\rm Aut}}
\newcommand{\arccot}{\mathop{\rm arccot}}
\newcommand{\capp}{{\mathop{\rm cap}}}
\newcommand{\rcapp}{{\mathop{\rm rcap}}}
\newcommand{\diam}{\mathop{\rm diam}}
\newcommand{\divv}{\mathop{\rm div}}
\newcommand{\codim}{\mathop{\rm codim}}
\newcommand{\RRe}{\mathop{\rm Re}}
\newcommand{\IIm}{\mathop{\rm Im}}
\newcommand{\Tr}{{\mathop{\rm Tr}}}
\newcommand{\Vol}{{\mathop{\rm Vol}}}
\newcommand{\card}{{\mathop{\rm card}}}
\newcommand{\supp}{\mathop{\rm supp}}
\newcommand{\sgn}{\mathop{\rm sgn}}
\newcommand{\essinf}{\mathop{\rm ess\,inf}}
\newcommand{\esssup}{\mathop{\rm ess\,sup}}
\newcommand{\Int}{\mathop{\rm Int}}
\newcommand{\Leibniz}{\mathop{\rm Leibniz}}
\newcommand{\lcm}{\mathop{\rm lcm}}
\newcommand{\loc}{{\rm loc}}

\newcommand{\mod}{\mathop{\rm mod}}
\newcommand{\spann}{\mathop{\rm span}}
\newcommand{\one}{1\hspace{-4.5pt}1}

\newcommand{\DWR}{}

\hyphenation{groups}
\hyphenation{unitary}

\newcommand{\cb}{{\cal B}}
\newcommand{\cc}{{\cal C}}
\newcommand{\cd}{{\cal D}}
\newcommand{\ce}{{\cal E}}
\newcommand{\cf}{{\cal F}}
\newcommand{\ch}{{\cal H}}
\newcommand{\ci}{{\cal I}}
\newcommand{\ck}{{\cal K}}
\newcommand{\cl}{{\cal L}}
\newcommand{\cm}{{\cal M}}
\newcommand{\cn}{{\cal N}}
\newcommand{\co}{{\cal O}}
\newcommand{\cs}{{\cal S}}
\newcommand{\ct}{{\cal T}}
\newcommand{\cx}{{\cal X}}
\newcommand{\cy}{{\cal Y}}
\newcommand{\cz}{{\cal Z}}

\newcommand{\wtozp}{W^{1,2}\raisebox{10pt}[0pt][0pt]{\makebox[0pt]{\hspace{-34pt}$\scriptstyle\circ$}}}
\newlength{\hightcharacter}
\newlength{\widthcharacter}
\newcommand{\covsup}[1]{\settowidth{\widthcharacter}{$#1$}\addtolength{\widthcharacter}{-0.15em}\settoheight{\hightcharacter}{$#1$}\addtolength{\hightcharacter}{0.1ex}#1\raisebox{\hightcharacter}[0pt][0pt]{\makebox[0pt]{\hspace{-\widthcharacter}$\scriptstyle\circ$}}}
\newcommand{\cov}[1]{\settowidth{\widthcharacter}{$#1$}\addtolength{\widthcharacter}{-0.15em}\settoheight{\hightcharacter}{$#1$}\addtolength{\hightcharacter}{0.1ex}#1\raisebox{\hightcharacter}{\makebox[0pt]{\hspace{-\widthcharacter}$\scriptstyle\circ$}}}
\newcommand{\scov}[1]{\settowidth{\widthcharacter}{$#1$}\addtolength{\widthcharacter}{-0.15em}\settoheight{\hightcharacter}{$#1$}\addtolength{\hightcharacter}{0.1ex}#1\raisebox{0.7\hightcharacter}{\makebox[0pt]{\hspace{-\widthcharacter}$\scriptstyle\circ$}}}

 \thispagestyle{empty}
  \begin{center}
% \vspace*{-1.0cm}
% {\Huge PRELIMINARY DRAFT 31/10/13}
\vspace*{1.5cm}

{\Large{\bf Gaussian bounds, strong ellipticity }}\\[3mm] 
{\Large{\bf and uniqueness criteria }}  \\[5mm]
\large Derek W. Robinson$^\dag$ \\[2mm]

\normalsize{November  2013}
\end{center}

\vspace{5mm}

\begin{center}
{\bf Abstract}
\end{center}

\begin{list}{}{\leftmargin=1.7cm \rightmargin=1.7cm \listparindent=15mm 
   \parsep=0pt}
   \item
Let $h$ be a quadratic form with domain $W_0^{1,2}(\Ri^d)$ given by
\[
h(\varphi)=\sum^d_{i,j=1}(\partial_i\varphi,c_{ij}\,\partial_j\varphi) 
\]
where $c_{ij}=c_{ji}$ are real-valued, locally bounded, measurable functions and 
 $C=(c_{ij})\geq 0 $. 
 If $C$ is strongly elliptic, i.e.\ if there exist $\lambda, \mu>0$ such that 
$\lambda\,I\geq C\geq \mu \,I>0$, then  $h$ is closable, the closure determines 
a positive self-adjoint operator $H$ on $L_2(\Ri^d)$ which generates a submarkovian semigroup $S$ with a  positive distributional kernel~$K$ and the kernel satisfies Gaussian upper and lower bounds.
Moreover, $S$ is conservative, i.e.\ $S_t\one=\one$ for all $t>0$.
Our aim is to examine converse statements.

First we establish that   $C$ is strongly elliptic if  and only if $h$ is closable, the  semigroup $S$ is conservative and $K$ satisfies Gaussian bounds.
Secondly, we prove that if the coefficients are such that a  Tikhonov growth condition is satisfied then $S$ is conservative. Thus in this case strong ellipticity of $C$ is equivalent to closability of $h$ together with Gaussian bounds on $K$.
Finally we consider coefficients $c_{ij}\in W^{1,\infty}_{\rm loc}(\Ri^d)$.
It follows that  $h$ is closable and a growth condition of the T\"acklind type is sufficient to establish the equivalence of  
strong ellipticity of $C$ and  Gaussian bounds on $K$.

\end{list}

\vfill

\noindent AMS Subject Classification: 35J15, 47B25, 47D07.

\vspace{0.5cm}

\noindent
\begin{tabular}{@{}cl@{\hspace{10mm}}cl}
$ {}^\dag\hspace{-5mm}$&   Mathematical Sciences Institute (CMA)    &  {} &{}\\
  &Australian National University& & {}\\
&Canberra, ACT 0200 && {} \\
  & Australia && {} \\
  &derek.robinson@anu.edu.au%& {}
 & &{}\\
\end{tabular}

\newpage

\setcounter{page}{1}

\newpage

\section{Introduction}\label{S1}

It is nearly fifty years since Aronson  \cite{Aro} established that the fundamental solutions of 
parabolic diffusion equations  are bounded above and below by Gaussian functions.
The essential  assumption in Aronson's argument is  strong ellipticity
of the  elliptic part of  the parabolic equation.
 Our aim is to examine converse statements, i.e.\  to analyze conditions which  ensure that Gaussian bounds imply strong ellipticity.
 We consider second-order operators in divergence form  acting on $\Ri^d$ determined by a $d\times d$-matrix
 $C=(c_{ij}) $ of  coefficients $c_{ij}$.
 We assume throughout that $c_{ij}=c_{ji}$ are real-valued, locally bounded, measurable functions
 and that $C(x)\geq 0$ for almost all $x\in \Ri^d$.
 Our principal results also  require bounds, either explicit or implicit, on the growth of the matrix norm $\|C(x)\|$
 as $x\to\infty$.
 But no growth conditions are necessary for our initial formulation.
 
We begin by introducing the Markovian form $h$ by the definition
  \begin{equation}
 \left\{ \begin{array}{ll}
\hspace{5mm}D(h)&\!\!=\;W_0^{1,2}(\Ri^d)\\[5pt]
\hspace{5mm}h(\varphi)&\!\!=\;\sum^d_{i,j=1}(\partial_i\varphi,c_{ij}\,\partial_j\varphi)
 \end{array}\right.\label{ehrt1.2}
\end{equation}
The form is not necessarily closed nor even closable.
If, however, $h$ is closable then the closure ${\overline h}$  is automatically a local Dirichlet form.
This form determines 
a positive self-adjoint operator $H$ on $L_2(\Ri^d)$, formally identifiable as   $-\sum^d_{i,j=1}\partial_ic_{ij}\partial_j$,
which generates a submarkovian semigroup $S$ with a  positive distributional kernel~$K$.
Specifically  $S$ extends from $L_2(\Ri^d)\cap L_p(\Ri^d)$  to a  semigroup on each of the $L_p$-spaces
which is strongly continuous if $p\in[1,\infty\rangle$ and weakly$^*$ continuous if $p=\infty$.
Moreover,  if $0\leq \varphi\leq \one$ then $0\leq S_t\varphi\leq \one$.

The coefficient matrix $C$ and the corresponding form $h$ are defined to be strongly elliptic if there exist $\lambda, \mu>0$ such that 
\begin{equation}
\lambda\,I\geq C(x)\geq \mu \,I>0\;,
\label{ehrt1.1}
\end{equation}
for almost all $x\in \Ri^d$, in the sense of matrix ordering.
Strong ellipticity immediately implies that $h$ is closable and $D(\overline h)=W^{1,2}(\Ri^d)$.
Moreover, Aronson's arguments \cite{Aro}   establish that there are 
$a,a',b,b'>0$ such that
\begin{equation}
a\,G_{b;t}(x-y)\geq K_t(x\,;y)\geq a'\,G_{b'\!;t}(x-y)
\label{ehrt1.3}
\end{equation}
for  almost  all $x,y\in \Ri^d$ and all $t>0$ where $G_{b;t}(x)=t^{-d/2}\,e^{-b\,|x|^2t^{-1}}$.
In recent  years it has become increasingly clear that these Gaussian bounds encapsulate a great deal
 of information concerning the solutions of the corresponding diffusion equations or the  semigroup kernels.
For example, Fabes and Stroock \cite{FaS} demonstrated that the Gaussian bounds were sufficient to derive
the Nash--De Giorgi \cite{Nash} \cite{DG} results on the local H\"older continuity of the solutions.
The bounds can also be used to deduce that the semigroup $S$ is conservative, i.e.\ $S_t\one=\one$ for all $t>0$,
although this property was studied well before Aronson's work (see, for example, \cite{Gaf}).

The following statement is the simplest version of our results.

\begin{thm}\label{thrt1.1}
Assume $h$ is the Markovian form  with $L_{\infty,\rm loc}$-coefficients defined by $(\ref{ehrt1.2})$.
Then the  following conditions are equivalent:
\begin{tabel}
\item\label{thrt1.1-1}
the matrix of coefficients $C$ is strongly elliptic,
\item\label{thrt1.1-2}
the form $h$ is closable, the  associated semigroup $S$ is conservative and the 
semigroup kernel $K$ satisfies the Gaussian bounds $(\ref{ehrt1.3})$.
\end{tabel}
\end{thm}

A more detailed version of this result, Theorem~\ref{thrt2.1}, in which the assumption that $h$ is closable is circumvented
by introduction of the relaxation of  $h$ and in which the Gaussians are replaced by more general functions will be established in Section~\ref{S2}.
The  new element of the theorem, the implication  \ref{thrt1.1-2}$\Rightarrow$\ref{thrt1.1-1},
  strengthens  an earlier result of \cite{ERZ1} for operators with bounded  coefficients.
It might appear surprising that 
this implication   only appears to require local boundedness of the coefficients.
In fact the conservation condition for $S$  places an implicit  restriction on the growth of the coefficients at infinity.
This will be discussed in Sections~\ref{S3} and~\ref{S4}.

In Section~\ref{S3} we establish that the conservation property is satisfied if the  coefficient growth is limited by a version of  the Tikhonov condition \cite{Tik} for uniqueness of solutions of parabolic equations.
Then in Section~\ref{S4} we demonstrate that  for locally Lipschitz coefficients a stronger growth is possible.
The limit  on growth is dictated by a version of the  T\"acklind  uniqueness condition  \cite{Tac}.
These  conditions,  in combination with Theorems~\ref{thrt1.1} or~\ref{thrt2.1}, give alternative characterizations of strong ellipticity (see Theorems~\ref{thrt3.1} and~\ref{thrt4.1}).

\section{Strong ellipticity}\label{S2}

In this section we discuss the proof of Theorem~\ref{thrt1.1} and its variants.
Since the implication \ref{thrt1.1-1}$\Rightarrow$\ref{thrt1.1-2} is classical we concentrate on the converse.

First we reformulate  the theorem in terms of the relaxation $h_0$ of $h$.
The relaxation occurs in convergence theory \cite{ET}
\cite{DalM} and  is variously referred to as the relaxed form or  lower semi-continuous regularization.
 In our earlier papers \cite{ERSZ1} \cite{ERSZ2} \cite{ERZ1} we also used the term viscosity closure as one possible definition of the relaxation is analogous to the viscosity approximation method of partial differential equations (see \cite{ERSZ1}, Section~2).
 The following definition is, however,  most suited to our immediate purposes.

The relaxation $h_0$ of the form $h$ is the  largest positive  closed quadratic form which is majorized by $h$, i.e.\
it is the largest closed quadratic form with $D(h)\subseteq D(h_0)$ satisfying $h_0(\varphi)\leq h(\varphi)$ for all $\varphi\in D(h)$.
In particular  if $h$ is closable then $h_0$ is equal to the closure~$\overline h$.
Note that if $h$ is not closable then $h_0$ is not an extension of $h$.

Throughout the following we let $H_0$ denote the positive  self-adjoint operator corresponding to $h_0$.
 Further, 
$S^{(0)}$ denotes the submarkovian semigroup generated by $H_0$ and~$K^{(0)}$ 
 the distributional kernel of $S^{(0)}$.

 Theorem~\ref{thrt1.1} now has the following gneralization.

\begin{thm}\label{thrt2.1}
Assume $h$ is the Markovian form  with $L_{\infty,\rm loc}$-coefficients defined by $(\ref{ehrt1.2})$.
Let $h_0$  denote the relaxation of $h$, $S^{(0)}$  the associated semigroup and $K^{(0)}$ the semigroup kernel.
Then the following conditions are equivalent:
\begin{tabel}
\item\label{thrt2.1-1}
the matrix of coefficients $C$ is strongly elliptic,
\item\label{thrt2.1-2}
 the  semigroup $S^{(0)}$ is conservative and there exist two positive, non-zero, bounded functions $\sigma, \rho$ 
 with $\int_0^\infty dr \, r^{(d+1)/2}\sigma(r)< \infty$ and $ \int_0^\infty dr \, r^{(d+1)/2}\rho(r)>0$ such that 
 \begin{equation}
t^{-d/2}\sigma(|x-y|^2t^{-1})\geq  K^{(0)}_t(x\,;y)\geq t^{-d/2}\rho(|x-y|^2t^{-1})
 \label{ehrt2.10}
 \end{equation}
 for almost all $x,y\in\Ri^d$ and for all small $t>0$.
\end{tabel}
\end{thm}
\proof\ 
First note that \ref{thrt1.1-1}$\Rightarrow$\ref{thrt1.1-2} in  Theorem~\ref{thrt1.1} is the classic result
and Condition~\ref{thrt1.1-2} in  Theorem~\ref{thrt1.1} obviously implies Condition~\ref{thrt2.1-2} in  Theorem~\ref{thrt2.1}.
Therefore if we prove that the latter condition implies strong ellipticity of $C$ then this establishes both theorems.

We begin by deducing the lower bound of the strong ellipticity condition (\ref{ehrt1.1}) from the lower kernel bound.

\begin{prop}\label{phrt2.1} 
Assume the lower kernel bound  of $(\ref{ehrt2.10})$.
Then there exists a $\mu>0$ such that $C\geq \mu I$ almost everywhere.
\end{prop}
 \proof\
 The proof is a repetition of the argument given in \cite{ERZ1}.
 It  does not assume that $S^{(0)}$ is conservative.
 Since a similar argument is needed for the subsequent discussion of upper bounds we recall the details.
 
 The starting point is the observation that 
 \[
 h_0(\varphi)\geq t^{-1}(\varphi, (I-S^{(0)}_t)\varphi) 
 \]
 for all $\varphi\in D( h_0)$  by spectral theory.
But $S^{(0)}$ is submarkovian so $0\leq S^{(0)}_t\one\leq \one$ for all $t>0$. 
Therefore 
\begin{eqnarray*}
h_0(\varphi)
& \geq & (2t)^{-1}\Big( (S^{(0)}_t\one,|\varphi|^2)+(|\varphi|^2,S^{(0)}_t\one)
    - (\varphi, S^{(0)}_t\varphi) - (S^{(0)}_t\varphi,\varphi) \Big) 
\end{eqnarray*}
for all $\varphi\in D(h_0)$ and $t>0$ where $(\,\cdot\,,\,\cdot\,)$ denotes the duality pairing between $L_p$ and $L_q$.

Next let $\varphi,\chi_n \in C_c^\infty(\Ri^d)$
where the $\chi_n$ are an increasing family with  $0\leq \chi_n\leq1$,  $\chi_n=1$ on the support of $\varphi$
and  such  that $\chi_n\to1$  as $n\to\infty$.
Then since $S^{(0)}$ is positive $S^{(0)}_t\chi_n\leq S^{(0)}_t\one$ and
\begin{eqnarray*}
h_0(\varphi)
& \geq & (2t)^{-1}\Big( (S^{(0)}_t\chi_n,\chi_n|\varphi|^2)+(\chi_n|\varphi|^2,S^{(0)}_t\chi_n)
    - (\chi_n\varphi, S^{(0)}_t\chi_n\varphi) - (S^{(0)}_t\chi_n\varphi,\chi_n\varphi) \Big) \\[5pt]
&=& (2t)^{-1} \int_{\Ri^d }dx\int_{\Ri^d} dy \,
       K^{(0)}_t(x\,;y)\,\chi_n(x) \, \chi_n(y) \, |\varphi(x)-\varphi(y)|^2
\end{eqnarray*}
for all $t>0$.
Hence it follows from the lower bound in (\ref{ehrt2.10}) that
\[
h_0(\varphi)
\geq (2t)^{-1}\int_{\Ri^d}dx\int_{\Ri^d}dy\,t^{-d/2}\,\rho(|x-y|^2t^{-1})\, 
\chi_n(x) \, \chi_n(y) \, |\varphi(x)-\varphi(y)|^2 
\;\;\; .  
\]
Therefore, in the limit $n\to\infty$,  one deduces   that 
\[
h_0(\varphi)
\geq (2t)^{-1}\int_{\Ri^d}dx\int_{\Ri^d}dy\,t^{-d/2}\,\rho(|x-y|^2t^{-1})\, 
|\varphi(x)-\varphi(y)|^2 
\]
for all $\varphi\in C_c^\infty(\Ri^d)$ and $t\in\langle0,1]$.
If ${\widetilde\varphi}$ denotes the Fourier transform of $\varphi$ one then calculates,
with a change of variables,  that
\begin{eqnarray}
h_0(\varphi)
&\geq& 
t^{-1}\int_{\Ri^d}dx\,t^{-d/2}\,\rho(|x|^2t^{-1})\int_{\Ri^d}d\xi\, |{\widetilde\varphi}(\xi)|^2\,
(1-\cos\xi.x) \nonumber \\[5pt]
&=& 
t^{-1}\int_{\Ri^d}dx\,\rho(|x|^2)\int_{\Ri^d}d\xi\, |{\widetilde\varphi}(\xi)|^2 \,
(1-\cos t^{1/2}\xi.x) \nonumber  \\[5pt]
&=&2\int_{\Ri^d}d\xi\, |{\widetilde\varphi}(\xi)|^2\,\int_{\Ri^d}dx\,\rho(|x|^2)\,
t^{-1} \sin^2 (2^{-1} t^{1/2} \xi.x) \nonumber 
\end{eqnarray}
for all $\varphi\in C_c^\infty(\Ri^d)$ and $t\in\langle0,1]$.
Thus  in the limit $t\to0$ one has
\[
h_0(\varphi)
\geq 2^{-1}\int_{\Ri^d}d\xi\, |{\widetilde\varphi}(\xi)|^2\,\int_{\Ri^d}dx\,\rho(|x|^2)\, (\xi.x)^2
=\mu\int_{\Ri^d} d\xi\, |{\widetilde\varphi}(\xi)|^2\, |\xi|^2= \mu\,\|\nabla\varphi\|_2^2
\]
for all $\varphi\in C_c^\infty(\Ri^d)$
with $\mu=(2\,d)^{-1}\int_{\Ri^d}dx\,\rho(|x|^2)|x|^2$.
Then by another  change of variables one has $\mu\sim \int^\infty_0dr\,r^{(d+1)/2}\rho(r)>0$.
Since $h\geq h_0$ by the definition of the relaxation it also follows that  $h(\varphi)\geq \mu\,\|\nabla\varphi\|_2^2
$ for all $\varphi\in C_c^\infty(\Ri^d)$.

 Now for each $\varphi\in C_c^\infty(\Ri^d)$, $\xi\in \Ri^d$ and $k\in\Ri$ introduce the $C_c^\infty$-functions $\varphi_{c,k}, \varphi_{s,k}$ by $\varphi_{c,k}(x)=\cos(k\,x.\xi)\varphi(x)$ and $\varphi_{s,k}(x)=\sin(k\,x.\xi)\varphi(x)$. 
Then
\begin{eqnarray*}
\int_{\Ri^d}dx\,(\xi,C(x)\xi)\,|\varphi(x)|^2&=&\lim_{k\to\infty}k^{-2}(h(\varphi_{c,k})+ h(\varphi_{s,k}))\\[5pt]
&\geq&
\mu \lim_{k\to\infty} k^{-2}(\|\nabla\varphi_{c,k}\|_2^2+ \|\nabla\varphi_{s,k}\|_2^2)=\mu \int_{\Ri^d}dx\,|\xi|^2\,|\varphi(x)|^2
\;.
\end{eqnarray*}
Therefore $C\geq \mu I$ almost everywhere.
\hfill$\Box$

\bigskip

There is a similar implication for the upper ellipticity bound but this requires  all the elements of Condition~\ref{thrt2.1-2}
of Theorem~\ref{thrt2.1}.

 \begin{prop}\label{phrt2.2}
 Assume $S^{(0)}$ is conservative and that  the  kernel bounds   $(\ref{ehrt2.10})$ are satisfied.
Then there exists a $\lambda<\infty$ such that $C\leq \lambda I$ almost everywhere.
 \end{prop}
 \proof\  First,  since  $S^{(0)}$ is conservative one has
\[
t^{-1}(\varphi, (I-S^{(0)}_t)\varphi) =(2\,t)^{-1}\Big( (S^{(0)}_t\one,|\varphi|^2)+(|\varphi|^2,S^{(0)}_t\one)
    - (\varphi, S^{(0)}_t\varphi) - (S^{(0)}_t\varphi,\varphi) \Big) 
\]
for all $\varphi\in L_2(\Ri^d)$ and $t>0$.

Secondly,  let $\varphi, \chi_n\in C_c^\infty(\Ri^d)$ where the $\chi_n$ are  the approximation to the identity used in the proof of the preceding proposition.
Then 
\begin{eqnarray*}
t^{-1}(\varphi, (I-S^{(0)}_t)\varphi) &=&
\lim_{n\to\infty}(2\,t)^{-1}\Big( (\chi_n,S^{(0)}_t\chi_n|\varphi|^2)+(S^{(0)}_t\chi_n|\varphi|^2,\chi_n)\\[3pt]
&&\hspace{6cm}{}
    - (\chi_n\varphi, S^{(0)}_t\chi_n\varphi) - (S^{(0)}_t\chi_n\varphi,\chi_n\varphi) \Big) \\[5pt]
&=& \lim_{n\to\infty}(2\,t)^{-1} \int_{\Ri^d }dx\int_{\Ri^d} dy \,
       K^{(0)}_t(x\,;y)\,\chi_n(x) \, \chi_n(y) \, |\varphi(x)-\varphi(y)|^2
   \;.    
   \end{eqnarray*}
 But using the upper  bound in (\ref{ehrt2.10}) 
 one deduces that 
\[
t^{-1}(\varphi, (I-S^{(0)}_t)\varphi) \leq 
a\,(2\,t)^{-1}  \int_{\Ri^d }dx\int_{\Ri^d} dy \,
t^{-d/2}\,\sigma(|x-y|^2t^{-1})
\,|\varphi(x)-\varphi(y)|^2
\]
for all $t\in\langle0,1]$.
Then, however, one can repeat the reasoning used in the proof of Proposition~\ref{phrt2.1} to conclude that
\[
t^{-1}(\varphi, (I-S^{(0)}_t)\varphi) \leq 
2\,a\,\int_{\Ri^d}d\xi\, |{\widetilde\varphi}(\xi)|^2\,\int_{\Ri^d}dx\,\sigma(|x|^2)\,
t^{-1} \sin^2 (2^{-1} t^{1/2} \xi.x) 
\;.
\]
Therefore
\[
h_0(\varphi)=\lim_{t\to0}t^{-1}(\varphi, (I-S^{(0)}_t)\varphi) 
\leq 
\lambda \int_{\Ri^d} d\xi\, |{\widetilde\varphi}(\xi)|^2\, |\xi|^2= \lambda\,\|\nabla\varphi\|_2^2
\]
for all $\varphi\in C_c^\infty(\Ri^d)$ with $\lambda=a\, (2\,d)^{-1}\int_{\Ri^d}dx\,\sigma(|x|^2)|x|^2\sim \int^\infty_0 dr\,r^{(d+1)/2}\sigma(r)<\infty$.

Thirdly, it follows from  the lower  bound in (\ref{ehrt2.10}), by Proposition~\ref{phrt2.1},  that the lower ellipticity bound $C\geq \mu I>0$ is valid.
 But this is sufficient to deduce that $h$ is closable (see \cite{FOT}, Section~3.1, or \cite{MR}, Section~II.2a).
 Therefore $h_0=\overline h$.
 Hence the foregoing bounds give $h(\varphi)\leq \lambda\,\|\nabla\varphi\|_2^2$ for all $\varphi\in C_c^\infty(\Ri^d)$.

Finally it follows by the  argument used to complete the proof of Proposition~\ref{phrt2.1} that the upper ellipticity bound $C\leq \lambda I$ is valid.
\hfill$\Box$

\bigskip

The proof of the implication \ref{thrt2.1-2}$\Rightarrow$\ref{thrt2.1-1} in Theorem~\ref{thrt2.1}  is now a corollary of Propositions~\ref{phrt2.1}  and \ref{phrt2.2}.
This completes the proofs of both Theorems~\ref{thrt2.1} and \ref{thrt1.1}.
\hfill$\Box$

\section{The conservation property}\label{S3}

The characterizations of strong ellipticity given by Theorems~\ref{thrt1.1} and \ref{thrt2.1} depend on the relaxation semigroup
$S^{(0)}$ being conservative.
In this section we establish that this property follows from a growth bound on the coefficients $c_{ij}$.
Consequently one obtains more direct characterizations  of strong ellipticity.

The conservation property has been widely studied for second-order operators with smooth coefficients acting on Riemannian  manifolds (see, for example, \cite{Gri7} for background information).
Our discussion follows in part an argument originating with Gaffney \cite{Gaf} which was subsequently clarified and extended by Davies \cite{Dav12}.
(The latter article also contains many references to the earlier literature.)
Since we are dealing with the relaxed forms associated to degenerate operators with measurable coefficients several additional problems arise which complicate the arguments.
Therefore it is  convenient to adopt a different characterization of the relaxation $h_0$ of the form $h$.

Let  $h_\varepsilon$ denote the form with coefficients $C_\varepsilon=C+\varepsilon I$
and with   $\varepsilon>0$.
Since the $C_\varepsilon$ satisfy the lower ellipticity bound $C_\varepsilon\geq \varepsilon I>0$ the forms $h_\varepsilon$ are closable.
Moreover, their closures $\overline h_\varepsilon$ form a monotonically decreasing family of Dirichlet forms.
Then it follows from a result of Kato, \cite{Kat1} Theorem VIII.3.11, that the positive self-adjoint operators $H_\varepsilon$ corresponding to the $\overline h_\varepsilon$ converge in the strong resolvent sense to a positive self-adjoint operator $H_0$.
The closed quadratic form $h_0$ corresponding to $H_0$, i.e.\ the form with domain $D(h_0)=D(H_0^{1/2})$ given by
$h_0(\varphi)=\|H_0^{1/2}\varphi\|_2^2$, is the relaxation of $h$.
(For the equivalence of this definition with that given in the previous section see \cite{bSim4} \cite{bSim5}.)
The proof that $S^{(0)}$ is conservative can be achieved 
in two steps following a strategy  used in \cite{RSi2a} \cite{RSi2}  and \cite{OuR}.
First one proves that the submarkovian semigroups $S^{(\varepsilon)}$ associated with the Dirichlet forms $\overline h_{\varepsilon}$ are 
conservative.
This step is achieved by Davies' arguments \cite{Dav12}.
Secondly,  one establishes that the $S^{(\varepsilon)}$  are $L_1$-convergent to $S^{(0)}$ as $\varepsilon\to 0$.
Then  by taking adjoints one deduces that the $S^{(\varepsilon)}$ are weakly$^*$ convergent to $S^{(0)}$ on $L_\infty(\Ri^d)$ as $\varepsilon\to 0$.
Therefore $S^{(0)}$ is conservative by the first step.

The conservation property for $S^{(0)}$ is dependent on a growth condition  on the coefficients.
Let  $\nu(s)=\esssup_{|x|\leq s}\|C(x)\|$,
\[
\rho(s)=\int^s_0dt\,(1+\nu(t))^{-1/2}
\]
and  define the  balls $B_\rho(r)$ by
\[
B_\rho(r)=\{x\in\Ri^d: \rho(|x|)<r\}
\;.
\]
Then the growth condition  is the requirement that there are $a,b>0$ such that
\begin{equation}
 |B_\rho(r)|\leq a\,e^{b\,r^2}
\label{ehrt3.1}
\end{equation}
for all large $r\geq 1$ where $|B|$ denotes the Lebesgue measure of the set $B$.

Condition (\ref{ehrt3.1}) is the direct analogue of the bound defining the Tykhonov class of functions
\cite{Tik}  used in the discussion of uniqueness of solutions of the Cauchy equation (see \cite{Gri7}, Section~11.4).
It  is satisfied if $\|C(x)\|\leq c\, (1+|x|)^2\log(2+|x|)$ and this  is essentially  the maximal allowed  growth.

 Note that (\ref{ehrt3.1})  automatically implies the weaker growth condition 
\begin{equation}
\lim_{r\to \infty}\rho(r)=\infty
\label{ehrt3.2}
\end{equation}
 because if (\ref{ehrt3.2})  is false then $\rho$ is bounded and $B_\rho(r)=\infty$ for all large~$r$.

\begin{thm}\label{thrt3.1}
Assume $h$ is the Markovian form  with $L_{\infty,\rm loc}$-coefficients defined by $(\ref{ehrt1.2})$
 and  let $h_0$  denote the relaxation of $h$.
 If the  Tikhonov growth condition $(\ref{ehrt3.1})$ is satisfied then the semigroup $S^{(0)}$ 
 associated with $h_0$ is conservative.
\end{thm}
\proof\
First we consider the approximating forms $h_\varepsilon$ and the corresponding semigroups $S^{(\varepsilon)}$ with $\varepsilon\in\langle0,1]$.

\begin{lemma}\label{lhrt3.1}
Assume the growth condition $(\ref{ehrt3.2})$.
Then $D(\overline h_\varepsilon)=D_1$ where
\[
D_1=\{\varphi\in W^{1,2}_{\rm loc}(\Ri^d):\Gamma(\varphi)+|\nabla\varphi|^2+\varphi^2\in L_1(\Ri^d)\}
\]
and 
\[
\overline h_\varepsilon(\varphi)=\int_{\Ri^d}\Gamma_\varepsilon(\varphi)
\]
for each $\varphi\in D_1$ and $\varepsilon\in\langle0,1]$ where $\Gamma_\varepsilon(\varphi)=\Gamma(\varphi)+\varepsilon\,|\nabla\varphi|^2$.
\end{lemma}
\proof\
The forms $\varphi\in D_1\mapsto k_\varepsilon (\varphi)=\int_{\Ri^d}\Gamma_\varepsilon(\varphi)$
are closed (again see \cite{FOT}, Section~3.1, or \cite{MR}, Section~II.2a).
Moreover, $k_\varepsilon\supseteq h_\varepsilon$.
Therefore to deduce  that  $k_\varepsilon=\overline h_\varepsilon$ it is necessary to prove that $W^{1,2}_0(\Ri^d)$ is a core of $k_\varepsilon$.

Let $\eta\in C_c^\infty(\Ri)$ satisfy $0\leq \eta\leq1$, $\eta(s)=1$ if $s\in[0,1\rangle$,  $\eta(s)=0$ if $s\geq 2$ and $|\eta'|\leq 2$. 
Then define $\eta_R$ for $R>0$ by $x\in\Ri^d\mapsto \eta_R(x)=\eta(R^{-1}\rho(|x|))$.
It follows  from the growth assumption on  $\rho$ that the $\eta_R$ have support in the closure of $B_\rho(2R)$ and $\eta_R(x)\to1$  pointwise as $R\to\infty$.
Moreover, 
\[
0\leq \Gamma_\varepsilon(\eta_R)\leq \sup_{x\in\Ri^d}(1+\nu(|x|))\,|(\nabla\eta_R(x)|^2\leq 4\,R^{-2}\sup_{x\in\Ri^d}
(1+\nu(|x|))\,|(\nabla\rho)(|x|)|^2\leq 4\,R^{-2}
\]
for all $\varepsilon\in\langle0,1]$.
Hence $\Gamma_\varepsilon(\eta_R)(x)\to0$ as $R\to\infty$.

Now if  $\varphi\in D_1\cap L_\infty(\Ri^d)$ and  $\varphi_R=\eta_R\varphi$ then
 $\varphi_R\in W^{1,2}_0(\Ri^d)$.
 But
\[
\Gamma_1(\varphi_R)+\varphi_R^2\leq 2\,\Gamma_1(\eta_R)\,\varphi^2+2\,\Gamma_1(\varphi)+\varphi^2\in L_1(\Ri^d)
\;.
\]
Therefore $\varphi_R\in D_1$.
But $\|\varphi_R-\varphi\|_2^2\to0$ as $R\to\infty$
 by the dominated convergence theorem.
In addition
\begin{eqnarray*}
 k_\varepsilon(\varphi_R-\varphi)&\leq& 2\int_{\Ri^d}\Gamma_\varepsilon(\eta_R)\,\varphi^2
+2\int_{\Ri^d}(\eta_R-\one)^2\,\Gamma_\varepsilon(\varphi)
\end{eqnarray*}
and both terms on the right converge to zero as $n\to\infty$ by dominated convergence.
Hence $\varphi_R\to\varphi$ in the $D( k_\varepsilon)$-graph norm.
Consequently 
$D( k_\varepsilon)\cap L_\infty(\Ri^d)\subseteq D(\overline h_\varepsilon)\cap L_\infty(\Ri^d)$.
But since both the closed  forms are Dirichlet forms this suffices to deduce  their equality.
 \hfill$\Box$

\bigskip

The identification of the domains of the forms $\overline h_\varepsilon$ is the starting point for application of the Davies--Gaffney arguments to the semigroups $S^{(\varepsilon)}$.
The method utilizes a family of bounded multiplication operators $\tau\in\Ri\mapsto U_\tau\colon U_\tau\varphi=e^{\tau\psi}\varphi$ 
where  $\psi\in W^{1,\infty}(\Ri^d)$   with $\Gamma(\psi)\leq 1$.
Clearly $U_\tau W^{1,2}_0(\Ri^d)\subseteq W^{1,2}_0(\Ri^d)$ but a simple estimate gives
\[
\Gamma(U_\tau\varphi)\leq 2\,U_{2\tau}\!\left(\Gamma(\varphi)+\tau^2\varphi^2\right)
\]
for all $\varphi\in W^{1,2}_0(\Ri^d)$.
Therefore by continuity $U_\tau D_1\subseteq D_1$ and the last estimate extends to all $\varphi\in D_1$.
This allows one to make estimates on the approximating semigroups $S^{(\varepsilon)}$.
The following lemma is the analogue of Lemma~1 in \cite{Dav12}.

\begin{lemma}\label{lhrt3.2}
If  $\varphi\in W^{1,2}_0(\Ri^d)$ then 
\[
\|U_\tau S^{(\varepsilon)}_t\varphi\|_2\leq e^{\tau^2 t}\|U_\tau\varphi\|_2
\]
and
\[
\int^t_0ds\int_{\Ri^d}\,U_{2\tau}\,\Gamma_\varepsilon(S^{(\varepsilon)}_s\varphi)\leq 2\, e^{2\tau^2t}\|U_\tau\varphi\|_2^2
\]
for all $\tau\in\Ri$, $t>0$ and  $\varepsilon\in\langle0,1]$.
\end{lemma}
\proof\
The proof begins with  the identity
\[
{{d}\over{dt}}\|U_\tau S^{(\varepsilon)}_t\varphi\|_2^2=-2\,\overline h_\varepsilon(S^{(\varepsilon)}_t\varphi,U_{2\tau}S^{(\varepsilon)}_t\varphi)
=-2\int_{\Ri^d}\Gamma_\varepsilon(S^{(\varepsilon)}_t\varphi, U_{2\tau}S^{(\varepsilon)}_t\varphi)
\;.
\]
Then one computes that 
\begin{eqnarray*}
{{d}\over{dt}}\|U_\tau S^{(\varepsilon)}_t\varphi\|_2^2&=&-2\int_{\Ri^d}U_{2\tau}\Gamma_\varepsilon(S^{(\varepsilon)}_t\varphi)-4\,\tau \int_{\Ri^d}(U_{2\tau}S^{(\varepsilon)}_t\varphi)\,\Gamma_\varepsilon(S^{(\varepsilon)}_t\varphi,\psi)\\[5pt]
&\leq&-2\int_{\Ri^d}U_{2\tau}\Gamma_\varepsilon(S^{(\varepsilon)}_t\varphi)+2\,\Big(\delta^{-1}\int_{\Ri^d}U_{2\tau}\Gamma_\varepsilon(S^{(\varepsilon)}_t\varphi)
+\delta\,\tau^2\int_{\Ri^d}U_{2\tau}(S^{(\varepsilon)}_t\varphi)^2\Big)\\[5pt]
&=&-2\,(1-\delta^{-1})\int_{\Ri^d}U_{2\tau}\Gamma_\varepsilon(S^{(\varepsilon)}_t\varphi)+2\,\delta\,\tau^2\,\|U_\tau S^{(\varepsilon)}_t\varphi\|_2^2
\end{eqnarray*}
for all $\delta>0$.

First choose $\delta=1$. Then
\[
{{d}\over{dt}}\|U_\tau S^{(\varepsilon)}_t\varphi\|_2^2\leq 2\,\tau^2\,\|U_\tau S^{(\varepsilon)}_t\varphi\|_2^2
\]
and on integration this gives the first statement of the lemma.

Secondly choose $\delta=2$.
Then 
\[
{{d}\over{dt}}\|U_\tau S^{(\varepsilon)}_t\varphi\|_2^2+\int_{\Ri^d}U_{2\tau}\Gamma_\varepsilon(S^{(\varepsilon)}_t\varphi)\leq 4\,\tau^2\,\|U_\tau S^{(\varepsilon)}_t\varphi\|_2^2\leq 4\,\tau^2\,e^{2\tau^2 t}\|U_\tau \varphi\|_2^2
\]
and on integration one obtains
\[
\|U_\tau S^{(\varepsilon)}_t\varphi\|_2^2-\|U_\tau \varphi\|_2^2+\int^t_0ds\int_{\Ri^d}\,U_{2\tau}\Gamma_\varepsilon(S^{(\varepsilon)}_s\varphi)
\leq 2\,(e^{2\tau^2 t}-1)\,\|U_\tau \varphi\|_2^2
\;.
\]
Therefore
\[
\int^t_0ds\int_{\Ri^d}\,U_{2\tau}\Gamma_\varepsilon(S^{(\varepsilon)}_s\varphi)\leq 
(2\,e^{2\tau^2 t}-1)\,\|U_\tau \varphi\|_2^2\leq 2\,e^{2\tau^2 t}\,\|U_\tau \varphi\|_2^2
\]
and the second statement of the lemma is verified.
\hfill$\Box$

\begin{remarkn}\label{rhrt3.1}
The estimates of Lemma~\ref{lhrt3.2}  extend to all $\varphi\in L_2(\Ri^d)$.
This is evident for the first estimate.
But the second can be extended in two steps.
First, $W^{1,2}_0(\Ri^d)$ is a core of $D_1$,
equipped with the $D(\overline h_1)$-graph norm,  and the estimate extends by continuity.
Secondly, $S^{(\varepsilon)}$ is holomorphic on $L_2(\Ri^d)$ so  $S^{(\varepsilon)}_uL_2(\Ri^d)\subseteq D_1$ for all $u>0$. Therefore replacing $\varphi$ by $S^{(\varepsilon)}_u\varphi$ one deduces that
\[
\int^{t+u}_uds\int_{\Ri^d}\,U_{2\tau}\,\Gamma_\varepsilon(S^{(\varepsilon)}_s\varphi)\leq 2\, e^{2\tau^2t}\|U_\tau S^{(\varepsilon)}_u\varphi\|_2^2
\]
for all $\varphi\in L_2(\Ri^d)$.
The full estimate then follows in the limit $u\to0$.
\end{remarkn}

\noindent{\bf Proof of Theorem~\ref{thrt3.1}}$\;$ Let $\eta_R$ be the approximation to the identity used in the proof of 
Lemma~\ref{lhrt3.1}.
Fix $\varphi\in C_c^\infty(B_\rho(r))$.
 Then
\begin{eqnarray*}
|(\one, S^{(\varepsilon)}_t\varphi)-(\one, \varphi)|&\leq&
\lim_{R\to\infty}|(\eta_R, S^{(\varepsilon)}_t\varphi)-(\eta_R, \varphi)|\\[5pt]
&\leq&
\lim_{R\to\infty}\Big|\int^t_0ds\int_{\Ri^d}\Gamma_\varepsilon(\eta_R, S^{(\varepsilon)}_s\varphi)\Big|
\;.
\end{eqnarray*}
Therefore,  by the Cauchy--Schwarz inequality, one has
\begin{eqnarray*}
\Big|\int^t_0ds\int_{\Ri^d}\Gamma_\varepsilon(\eta_R, S^{(\varepsilon)}_s\varphi)\Big|
&\leq &\Big(\int^t_0ds\,\int_{\Ri^d}\,U_{2\tau}\,\Gamma_\varepsilon(S^{(\varepsilon)}_s\varphi)\Big)^{1/2}
\cdot\Big(\int^t_0ds\,\int_{\Ri^d}\,U_{-2\tau}\,\Gamma_\varepsilon(\eta_R)\Big)^{1/2}
\;.
\end{eqnarray*}
Next set $\psi(x)=\rho(x)\wedge R-r$ with $R>r$ and assume $\tau\geq0$.
 Then, by  the second estimate of Lemma~\ref{lhrt3.2},
 \[
\int^t_0ds\int_{\Ri^d}\,U_{2\tau}\,\Gamma_\varepsilon(S^{(\varepsilon)}_s\varphi)\leq 2\, e^{2\tau^2t}\|U_\tau\varphi\|_2^2\leq 
2\, e^{2\tau^2t}\|\varphi\|_2^2
\]
where the second bound follows because 
 $|U_\tau\varphi|\leq |\varphi|$ for $\varphi$ with support in $B_\rho(r)$.
But  $\Gamma_\varepsilon(\eta_R)$ has support in the annulus $A_\rho(R)=\{x\in\Ri^d: R\leq |x|\leq 2R\}$  and $\Gamma_\varepsilon(\eta_R)\leq 4\,R^{-2}$.
Therefore
\[
\int^t_0ds\,\int_{\Ri^d}\,U_{-2\tau}\,\Gamma_\varepsilon(\eta_R)\leq 4\,R^{-2}\,t\int_{A_\rho(R)}\,e^{-2\tau(\rho\wedge R-r)}
\leq 4\,t\, R^{-2}\,|B_\rho(2R)|\,e^{-2\tau (R-r)}
\;.
\]
Combining these estimates one deduces that 
\[
\Big|\int^t_0ds\int_{\Ri^d}\Gamma_\varepsilon(\eta_R, S^{(\varepsilon)}_s\varphi)\Big|\leq 2^{3/2}\,(R^2/t)^{-1/2}\,e^{\tau^2t-\tau (R-r)}\,|B_\rho(2R)|^{1/2}\|\varphi\|_2
\;
\]
Hence setting $\tau=(R-r)(2t)^{-1}$ one obtains
\[
\Big|\int^t_0ds\int_{\Ri^d}\Gamma_\varepsilon(\eta_R, S^{(\varepsilon)}_s\varphi)\Big|\leq 2\,(R^2/(4t))^{-1/2}\,e^{-(R-r)^2(4t)^{-1}}\,|B_\rho(2R)|^{1/2}\|\varphi\|_2
\;
\]
Then it follows from the Tikhonov condition (\ref{ehrt3.1}) that 
\[
\lim_{R\to\infty}\Big|\int^t_0ds\int_{\Ri^d}\Gamma_\varepsilon(\eta_R, S^{(\varepsilon)}_s\varphi)\Big|=0
\]
if $t<16 \,b$.
Hence, by a density argument, one deduces that $S^{(\varepsilon)}_t\one=\one $ for all small $t>0$ and then by the semigroup property for all $t>0$.
Thus the approximating semigroups $S^{(\varepsilon)}$ are conservative.

The second step of the proof consists of establishing that the $S^{(\varepsilon)}_t$ are $L_1$-convergent to $S^{(0)}_t$ as 
$\varepsilon\to0$.
Then it immediately follows that $S^{(0)}$ is conservative.
Therefore we now fix $\varphi\in L_1(B_\rho(r))\cap L_2(B_\rho(r))$ and note that 
\begin{eqnarray}
\|(S^{(\varepsilon_1)}_t-S^{(\varepsilon_2)}_t)\varphi\|_1&\leq &
\|\one_{B_\rho(R)}(S^{(\varepsilon_1)}_t-S^{(\varepsilon_2)}_t)\varphi\|_1
+\|\one_{B_\rho(R)^{\rm c}}S^{(\varepsilon_1)}_t\varphi\|_1+\|\one_{B_\rho(R)^{\rm c}}S^{(\varepsilon_2)}_t\varphi\|_1
\nonumber\\[5pt]
&\leq&|B_\rho(R)|^{1/2}\|(S^{(\varepsilon_1)}_t-S^{(\varepsilon_2)}_t)\varphi\|_2
+2\,\sup_{\varepsilon\leq1}\int_{B_\rho(R)^{\rm c}}|S^{(\varepsilon)}_t\varphi|
\label{ehrt4.11}
\end{eqnarray}
for all $R>r$ and $\varepsilon_1,\varepsilon_2\leq 1$.
Since $S^{(\varepsilon)}_t$ is $L_2$-convergent to $S^{(0)}_t$  as $\varepsilon\to 0$ for all $t>0$ it suffices for the $L_1$-convergence to prove that the last term in (\ref{ehrt4.11}) converges to zero as $R\to\infty$ for all $t$ in a finite interval $\langle0,t_0]$ uniformly for all $\varepsilon\leq1$.

Set
 $A_n=\{x: n\leq \rho(|x|)<n+1\}$.
Then 
\begin{eqnarray*}
\int_{A_n}|S^{(\varepsilon)}_t\varphi|&\leq& |B_\rho(n+1)|^{1/2}\,\bigg(\int_{A_n}|S^{(\varepsilon)}_t\varphi|^2\bigg)^{1/2}
\;.
\end{eqnarray*}
One may assume $R$ is an integer.
Let $S$ be a second integer with $S\geq R+1$.
Now with   $\psi=\rho\wedge S-r$ and $\tau>0$ one deduces from the first estimate of Lemma~\ref{lhrt3.2} that
\begin{eqnarray*}
\int_{A_n}|S^{(\varepsilon)}_t\varphi|^2&\leq&\int_{A_n}U_{-2\tau}|U_\tau S^{(\varepsilon)}_t\varphi|^2\\[5pt]
&\leq &e^{-2\tau(n-r)}\|U_\tau S^{(\varepsilon)}_t\varphi\|_2^2\leq e^{\tau^2t-2\tau(n-r)}\|\varphi\|_2^2
\end{eqnarray*}
for all $R\leq n\leq S$.
Combining these estimates with $\tau=(n-r)/(2t)$ and using the volume bound (\ref{ehrt3.1}) then gives the estimate
\[
\bigg(\int_{A_n}|S^{(\varepsilon)}_t\varphi|\bigg)^2\leq a\,e^{b(n+1)^2}e^{-(n-r)^2/(4t)}\|\varphi\|_2^2
\;.
\]
 Therefore one may choose $t_0>0$ and $a'>0$ such that 
 \[
 \int_{A_n}|S^{(\varepsilon)}_t\varphi|\leq a'\,e^{-b\,n^2}\,\|\varphi\|_2
 \]
 for all $t\in\langle0,t_0]$, $\varepsilon\in\langle0,1]$ and  $R\leq n\leq S$.
The value of $a'$ depends on $r$ but $t_0$ is independent of $r$.
Finally \begin{eqnarray*}
\int_{B_\rho(R)^{\rm c}}|S^{(\varepsilon)}_t\varphi|& \leq & \lim_{S\to\infty}\sum_{n=R}^S\int_{A_n}|S^{(\varepsilon)}_t\varphi|
    \leq a'\,\Big(\sum_{n\geq R}  e^{-b\,n^2}\Big)\,\|\varphi\|_2
     \end{eqnarray*}
 for  $t\in\langle0,t_0] $ uniformly for $\varepsilon\in\langle0,1]$.
 Thus the last term in (\ref{ehrt4.11}) converges to zero as $R\to\infty$ for all small $t>0$, all $\varphi
 \in L_1(B_\rho(r))\cap L_2(B_\rho(r))$ and all $r$.
 Consequently $S^{(\varepsilon)}_t$ is $L_1$-convergent to $S^{(0)}_t$  for all small $t$ and then  for all $t>0$ by the 
semigroup property.
Since the $S^{(\varepsilon)}$ are conservative it  follows that $S^{(0)}$ is conservative.
\hfill$\Box$

\bigskip

Theorem~\ref{thrt3.1} allows an alternative formulation of the earlier characterizations,
Theorems~\ref{thrt1.1} and \ref{thrt2.1},  of strong ellipticity.

\begin{thm}\label{thrt3.2}
Assume $h$ is the Markovian form  with $L_{\infty,\rm loc}$-coefficients defined by $(\ref{ehrt1.2})$
for which the Tikhonov growth condition $(\ref{ehrt3.1})$ is satisfied.
Then the following conditions are equivalent:
\begin{tabel}
\item\label{thrt3.2-1}
the matrix of coefficients $C$ is strongly elliptic,
\item\label{thrt3.2-2}
the form $h$ is closable and the semigroup kernel $K$ corresponding to the closure $\overline h$ of $h$ satisfies Gaussian bounds,
\item\label{3.2-3}
the semigroup kernel $K^{(0)}$ corresponding to the relaxation $h_0$ of $h$ satisfies the Gaussian bounds $(\ref{ehrt1.3})$.
\end{tabel}
\end{thm}
One can  also replace the Gaussian bounds by the more general bounds (\ref{ehrt2.10}) of Theorem~\ref{thrt2.1} without destroying the equivalences.

\section{Locally Lipschitz coefficients}\label{S4}

In this section we establish that the previous characterizations of strong ellipticity can be strengthened if the coefficients $c_{ij}$ are regular. 
In particular the Tikhonov growth condition can be replaced by a weaker condition of the T\"acklind type.
Throughout this  section we assume that $c_{ij}=c_{ji}\in W^{1,\infty}_{\rm loc}(\Ri^d)$ and that $C\geq 0$.
Then   $L=-\sum^d_{i,j=1}\partial_ic_{ij}\partial_j$ with domain $D(L)=C_c^\infty(\Ri^d)$ is a symmetric operator  on $L_2(\Ri^d)$ and 
the form $h$, defined by (\ref{ehrt1.2}), satisfies $h(\varphi)=(\varphi, L\varphi)$ for all $\varphi\in D(L)$.
It follows that $h$ is closable and  the operator $H$ corresponding to the closure $\overline h$
is the Friedrichs extension of $L$.
We continue to denote  by $S$ the  submarkovian semigroup generated by $H$  and  $K$  the semigroup kernel.
The principal advantage of the regularity of the coefficients is that one can exploit elliptic and parabolic regularity properties.
In particular one can deduce that the   submarkovian semigroup $S$ generated by $H$ is conservative under more general 
conditions than previously.

First the T\"acklind condition on the coefficients is defined as the growth restriction
\begin{equation}
\int_R^\infty dr\,r\,(\log |B_\rho(r)|)^{-1}=\infty
\label{ehrt4.1}
\end{equation}
for all large $R$ with the convention $1/\infty=0$.
The earlier Tikhonov condition (\ref{ehrt3.1}) is equivalent to the existence of a $c>0$ such that $r\,(\log|B_\rho(r)|)^{-1}\geq c\,r^{-1}$ for all large $r$.
Hence it is a stronger restriction than (\ref{ehrt4.1}).
The T\"acklind condition still implies, however, that $\rho(x)\to\infty$ as $|x|\to\infty$.

\begin{thm}\label{thrt4.1}
Assume $h$ is the Markovian form  with $L_{\infty,\rm loc}$-coefficients defined by $(\ref{ehrt1.2})$.
Assume the T\"acklind growth condition $(\ref{ehrt4.1})$ is satisfied.

Then $S$ is conservative and the following conditions are equivalent:
 \begin{tabel}
 \item\label{thrt4.1-1}
 the matrix of coefficients $C$ is  strongly  elliptic,
 \item\label{thrt4.1-2}
  the semigroup kernel $K$ satisfies the Gaussian bounds $(\ref{ehrt1.3})$.
\end{tabel}
\end{thm} 
Again we note that the Gaussian bounds in Condition~\ref{thrt4.1-2} can be replaced by the more general bounds 
(\ref{ehrt2.10}) of Theorem~\ref{thrt2.1} without destroying the equivalence.

\medskip

\noindent{\bf Proof}$\;$ The implication  \ref{thrt4.1-1}$\Rightarrow$\ref{thrt4.1-2} 
is the classic Aronson result.
Moreover, once one establishes that $S$ is conservative the converse implication follows from Theorem~\ref{thrt1.1}.

The proof of  the conservation property of $S$
 is of a quite different nature to that of Section~\ref{S3}.
It is based on various uniqueness criteria for locally strongly elliptic operators and uses results and techniques developed for the analysis of the Laplace-Beltrami operator on a Riemannian manifold.
In particular it uses  an old argument of Grigor'yan for the Laplace-Beltrami operator on a geodesically complete manifold (see \cite{Gri6}, Theorem~1, or \cite{Gri5}, Theorem~9.1)
as adapted in a recent paper \cite{Rob8}  to the analysis of  locally strongly elliptic operators on domains of $\Ri^d$.

There are two key uniqueness properties, $L_1$-uniqueness and Markov uniqueness.
First, the operator $L$, viewed as an operator on $L_1(\Ri^d)$, is defined to be $L_1$-unique if it has a unique $L_1$-extension
 which generates an $L_1$-continuous semigroup.
 It follows by a version of the Lumer--Phillips theorem due to Arendt (see \cite{Ebe} Theorem~A.1.2) that this is the case if and only if  the $L_1$-closure $\overline L^{\,\scriptscriptstyle 1}$ of $L$ is a semigroup generator.
In which case $\overline L^{\,\scriptscriptstyle 1}$  is  the generator of the submarkovian semigroup $S$ acting on $L_1(\Ri^d)$.
Secondly, the operator $L$, acting  on $L_2(\Ri^d)$,  is defined to be Markov unique if  it has a unique self-adjoint extension which generates a submarkovian semigroup.
This is clearly the case if and only if the Friedrichs extension $H$ is the unique submarkovian generator.

Next note that if the semigroup kernel $K$ satisfies the lower Gaussian bounds, or the  lower bound in 
(\ref{ehrt2.10}),  then the lower ellipticity bound $C\geq \mu I>0$ is valid  by Proposition~\ref{phrt2.1}.
But the coefficients $c_{ij}$ are locally bounded by assumption.
Therefore $C$ is locally strongly elliptic, i.e.\ for each compact subset $V$ there are $\mu_V,\lambda_V>0$ such that $\lambda_VI\geq C(x)\geq \mu_V I$ for all $x\in V$.
The following proposition summarizes some recent results relating uniqueness properties and semigroup conservation  for  locally strongly  elliptic operators.

\begin{prop}\label{phrt4.10}
Assume that $c_{ij}=c_{ji}\in W^{1,\infty}_{\rm loc}(\Ri^d)$ and  that
 the semigroup kernel satisfies the   lower bound $(\ref{ehrt2.10})$ of  Proposition~$\ref{phrt2.1}$.

\begin{tabel}
\item\label{phrt4.10-1}
If $L$ is $L_1$-unique then $S$ is conservative.
Conversely,  if $S$ is conservative
 and  $\rho(x)\to\infty$ as $|x|\to\infty$ then  $L$ is $L_1$-unique. 
\item\label{phrt4.10-2}
If $L$ is $L_1$-unique then $L$ is Markov unique.
Conversely, if $L$ is Markov unique and the T\"acklind condition $(\ref{ehrt4.1})$
is satisfied then $L$ is $L_1$-unique.
\end{tabel}
\end{prop} 
\proof\  \ref{phrt4.10-1}.
The first statement follows because  $(\one, L\varphi)=0$ for all $\varphi\in C_c^\infty(\Ri^d)$.
Therefore  $(\one, \overline L^{\,\scriptscriptstyle 1}\!\varphi)=0$ for all $\varphi\in D(\overline L^{\,\scriptscriptstyle 1})$. 
Since $L_1$-uniqueness implies that $\overline L^{\,\scriptscriptstyle 1}$ is the generator of $S$ on $L_1(\Ri^d)$ it follows  that $d(\one, S_t\varphi)/dt=0$ for all
$\varphi\in  D(\overline L^{\,\scriptscriptstyle 1})$.
Consequently $S$, acting on $L_\infty(\Ri^d)$,  is conservative.
This argument is independent of any growth restriction.
The converse statement is essentially due to Hasminski~\cite{Has} and Azencott~\cite{Aze} (see \cite{Dav14}, Theorem~2.2).
A detailed proof which only uses the growth property $\rho(x)\to\infty$ is given for a broader class of operators with lower order terms in \cite{OuR}, Theorem~3.6.
The argument relies on the local strong ellipticity of $C$. 
This implication is, however, not necessary in the proof of Theorem~\ref{thrt4.1}.

\smallskip

\ref{phrt4.10-2}. If $L$ has two distinct submarkovian extensions then these generate distinct submarkovian semigroups
and their $L_1$-generators are distinct. But both generators are extensions of $L$ acting on $L_1(\Ri^d)$.
So $L$ is not $L_1$-unique.
The converse is the difficult part of the proof.
It is a consequence of Grigor'yan's arguments (see \cite{Gri6} Theorem~2) which are described in detail in \cite{Gri7}, Section~11.4.
Grigor'yan deals with the Laplace-Beltrami operator on a geodesically complete manifold but his methods adapt
to Markov unique diffusion operators (see  \cite{Rob8}, Section~3 and in particular Corollary~3.3).
The proof is a quadratic form calculation involving the closed form $\overline h$ and it is essential for the calculation that the domain of the form is sufficiently large.
But in the current context this is assured by the Markov uniqueness.
In fact one has 
\begin{equation}
D( \overline h)=\{\varphi\in W^{1,2}_{\rm loc}(\Ri^d): \Gamma(\varphi)+\varphi^2\in L_1(\Ri^d)\}
\label{ehart4.2}
\end{equation}
(see below) and this  suffices  to carry through  Grigor'yan's arguments as explained in Section~3 of \cite{Rob8}.
The latter proof has to be slightly modified since we are  using the `distance function'  $\rho$ in place of the Riemannian distance and the balls $B_\rho(r)$ in place of the Riemannian balls $B(r)$.
Therefore it is appropriate  to redefine the Aronson auxiliary function 
in the proof by $\xi_t=\nu\,(\rho-r)^2(t-s)^{-1}$.
This new definition still satisfies the crucial bound (13) of \cite{Rob8} and the rest of the proof is essentially unchanged.
In fact it is simplified since $L$ is a pure second-order divergence-form operator and there are no lower order terms to estimate.
\hfill$\Box$

 \bigskip
\noindent{\bf Proof of Theorem~\ref{thrt4.1} continued} $\;$ It remains to prove that $S$ is conservative.
But, by Proposition~\ref{phrt4.10},  it suffices to prove that $L$ is Markov unique.
In particular Markov uniqueness and the T\"acklind condition imply that $L$ is $L_1$-unique,
by Proposition~\ref{phrt4.10}.\ref{phrt4.10-2},  and this in turn implies that $S$ is conservative, by Proposition~\ref{phrt4.10}.\ref{phrt4.10-1}.
But the Markov uniqueness  of $L$  is a corollary of \cite{RSi5}, Theorem~2.1, or \cite{Rob8}, Proposition~2.1.
The argument is as follows.

First define a form $h_N$ by specifying its domain $D(h_N)$ to be the set on the right hand side of (\ref{ehart4.2}) and then setting $ h_N(\varphi)=\int_{\Ri^d}\Gamma(\varphi)$ for all $\varphi\in D(h_N)$.
The form is closed because of the local strong ellipticity of $C$ (see \cite{OuR} Proposition~2.1) and it is a Dirichlet form by 
standard estimates.
But if $k$ is a Dirichlet form extension of $h$ it follows that $h\subseteq k\subseteq h_N$.
This is the crucial result given by  \cite{RSi5}, Theorem~2.1, or \cite{Rob8}, Proposition~2.1.
Then, however, one has the following variation 
 of Lemma~\ref{lhrt3.1}.
\begin{lemma} \label{lhrt4.1} Assume that $\rho(x)\to\infty$ as $|x|\to\infty$.
Then $\overline h=h_N$ and $h$ is Markov unique.
\end{lemma}

The proof is a repetition of the argument used to prove Lemma~\ref{lhrt3.1}.

\medskip

Finally, since the T\"acklind condition implies the condition $\rho(x)\to\infty$ as $|x|\to\infty$ one concludes
from the lemma that $h$ is Markov unique.  Hence $S$ is conservative
 and the proof of the theorem is complete. \hfill$\Box$

\bigskip

The foregoing proof that the semigroup $S$ is conservative depends explicitly on two uniqueness properties, $L_1$-uniqueness and Markov uniqueness, but it also depends implicitly on a third such property.
Grigor'yan's original argument as detailed in Chapters~8 and~11 of \cite{Gri7} is based on Theorem~8.18 which establishes the equivalence of semigroup conservation (stochastic completeness), $L_1$-uniqueness (the absence of non-trivial $\alpha$-harmonic functions) and the uniqueness of bounded solutions of the Cauchy equation.
Theorem~11.8 then derives the conservation property with the aid of the T\"acklind condition
by deducing that the condition suffices to establish that the Cauchy equation does have a unique bounded solution.
These results depend on parabolic regularity properties which follow from the Lipschitz property of the coefficients.
In the earlier case of measurable coefficients considered in Sections~\ref{S2} and~\ref{S3} these properties are no longer available and it is unclear whether there is an equivalence between the conservation property and uniqueness of bounded solutions to a weaker form of the Cauchy equation.

\end{document}